\documentclass[a4paper,leqno]{article}

\usepackage{amsmath,amssymb,amscd,cite,authblk,multicol,color}

\allowdisplaybreaks

\begin{document}

\begin{center}

\textbf{\Large Projective (or spin) representations of finite groups. III}

\bigskip

{Satoe Yamanaka\footnote{Department of Liberal Studies, 
National Institute of Technology, Nara College, 22 Yata-cho, Yamatokoriyama, Nara 639-1080, JAPAN}, 
Tatsuya Tsurii\footnote{Tokyo University of Information Sciences, 
Onaridai 4-1, Wakaba-ku, Chiba-shi, Chiba, 265-8501, JAPAN},
Itsumi Mikami\footnote{Hirai Mathematics Institute}, and 
Takeshi Hirai\footnote{Kyoto University, Yoshida-honmachi, Sakyo-ku, Kyoto 606-8501, JAPAN}
\footnote{Hirai Mathematics Institute, 22-5 Nakazaichi-cho, Iwakura, Sakyo, Kyoto 606-0027, JAPAN}}

\bigskip

\end{center}

\begin{abstract}
In previous 
papers I and II under the same title, we  
proposed a practical method called 
{\it Efficient stairway up to the Sky}, 
and  apply it to some typical finite 
groups $G$, with Schur multiplier $M(G)$ containing prime 
number 3, 
 to construct explicitly their representation groups $R(G)$,  
and then, 
to construct a complete set of representatives  
of linear IRs of $R(G)$, which gives naturally, 
through sectional 
restrictions, 
 a complete 
set of
 representatives of spin IRs of $G$. 
  In the present paper, 
we are concerned mainly with group $G=G_{39}$ of order 27 in 
a list of  
Tahara's paper, with 
$M(G)={\mathbb Z}_3\times {\mathbb Z}_3$.   
In this case, to arrive up to {\it the Sky}, we have 
 two steps of one-step efficient central 
extensions. 
By the 1st step, we obtain a covering group of order 81, 
 and 
by the 2nd step we arrive to $R(G)$ of order 243. 
At the 1st step, to construct explicitly 
a complete set of representatives of 
IRs of the group, we apply Mackey's induced 
representations,
  and at the 2nd step, 
with a help of this result, we apply 
so-called classical method for 
semidirect product groups given by Hirai 
 and arrive to 
 a complete list of 
IRs of $R(G)$. 
 Then, using explicit realization 
of these IRs, we can  
 compute their 
characters (called spin characters).

\bigskip

\medskip

\noindent
{\bf 2020 Mathematics Subject Classification:} Primary 20C25; Secondary 20F05, 20E99.

\noindent
{\bf Key Words:} representation group of finite group, projective representation, spin representation, spin character.
\end{abstract}

\begin{multicols*}{2}

{\bf 1. Rappels. }  
Let $G$ be a finite group. 
Our practical 
method  \lq\lq{\it Efficient stairway up 
to the Sky}\,\rq\rq, starts with a presentation 
of $G$ by a pair of {\bf Set of generators} 
and {\bf Set of fundamental relations}, 
to construct $R(G)$. 
  A central extension 
\begin{gather}
\label{2024-01-22-1}
\quad 
1 \longrightarrow A \longrightarrow H
\stackrel{\alpha}{\longrightarrow} G\longrightarrow 1\quad ({\rm exact})
\end{gather}
of $G$ by an abelian group $A$ is called 
an {\it efficient covering} if 
 $A$ is contained in 
the commutator subgroup $[H,H]$ of $H$. 
Take a pair $x,y\in G$ commuting with each other and of the same 
order $d>1$. Then take a {\it one-step efficient 
covering} given by the 
commuting pair $[x,y]=1$, in such a way that,  $A=\langle z\rangle$ 
is a cyclic group generated by 
 a central element $z$ of order $d$, given under a homomorphism 
$\alpha$ as   
\begin{gather}
\label{2024-01-22-2}
[\xi,\eta]=z,\quad \xi\stackrel{\alpha}{\rightarrow}  x,\quad
\eta \stackrel{\alpha}{\rightarrow} y,\;\quad 
z\stackrel{\alpha}{\rightarrow}1.
\end{gather}

When  $M(G)$ is known, we may arrive up to $R(G)$ 
by repeating one-step efficient central extensions step by step because  
$R(G)$ is characterized \cite[\S 2]{HMTY} 
as an efficient central extension $H$ 
of $G$ such that
\begin{gather}
\label{2024-01-25-11}
M(G) \subset Z(H)\cap [H,H], \quad H/M(G)\cong G.
\end{gather}

In cases treating here, for 
every one-step efficient central extensions, 
the resulting central extension 
$H$ is expressed as semidirect product $H=U\rtimes W$, 
and 
 so, construction 
of its irreducible representations is realized by Mackey's 
induced 
 representations or by 
 the classical method 
 in \cite{THirai}. 
We have studied  
cases where $M(G)$ contains 
prime number 2. 
Present series, I in \cite{HMTY}, 
II in \cite{TYMH} and III here, is 
a trial to other cases of prime number $p>2$, for
$p=3$, and  
 their details will appear elsewhere, with   
 thanks to computor-aided calculations 
\cite{Tsur}, \cite{WP} etc.

{\bf 2. Structure of $G=G_{39}$ and coverings. }  
The group $G=G_{39}$ is presented  
as follows: 
\\
\quad{\bf Set of generators:} \hspace{4ex}$\{a,\; b\}$;
\\
\quad{\bf Set of fundamental relations:} 
\begin{gather}
\label{2024-01-23-1}
a^3=b^3=(ab)^3=(ab^{-1})^3=1.
\end{gather} 
Then, 
we can introduce another presentation  of $G$ 
as follows \cite[Lemma 5.2]{HMTY}: 
put $c:=[b,a]$ and 
\\
\quad{\bf Set of generators:}
\quad$
\{x_1=b,\;x_2=c,\;x_3=a\}\,; 
$
\\
\quad{\bf Set of fundamental relations:}
\begin{eqnarray}
\label{2024-01-23-2}
\left\{
\begin{array}{ll}
x_1^{\,3}=x_2^{\,3}=x_3^{\,3}=1,
\\
x_2\;\;{\rm central},\;\;\;
[x_1,x_3]=x_2\,,
\end{array}
\right.
\end{eqnarray}
and, with these structural data, 
 the group $G_{39}$ is identified with 
the group $G_{27}^{\;\,3}$ with GAPIdentity [27,3]. 
Its center is $X_2:=\langle x_2\rangle$. 
Taking the commuting pair $[x_1,x_2]=1$, we construct 
one-step efficient central extension 
$\widetilde{G}\stackrel{\delta}{\rightarrow}G$ relative to
 it as follows: 
\\
\quad{\bf Set of generators:} \quad 
$\{\xi_1,\xi_2,\xi_3,z_{12}\}$;
\\[.4ex]
\quad{\bf Covering map:}\quad  
$
\left\{
\begin{array}{ll}
\mbox{\rm $
\xi_i\stackrel{\delta}{\rightarrow}x_i\;(1 \leqslant i\leqslant 3),$}
\\
\mbox{\rm $
[\xi_1,\xi_2]=z_{12}\stackrel{\delta}{\rightarrow}
[x_1,x_2]=1;$}
\end{array}
\right.
$ 
\\[.6ex]
\quad
{\bf Set of fundamental relations:} 
\begin{gather}
\label{2023-11-11-11}
\left\{
\begin{array}{ll} 
\xi_1^{\;3}=\xi_2^{\;3}=
\xi_3^{\;3}=1, 
\\
z_{12}=[\xi_1,\xi_2]\;\;{\rm central},\;\;z_{12}^{\;\;3}=1,
\\
\xi_2=[\xi_1, \xi_3],\quad[\xi_2,\xi_3]=1.
\end{array}
\right.
\end{gather}
We can summarize one of main results in 
\cite[\S  5]{HMTY} as

{\bf Theorem 2.1 (cf. \cite[Theorem 5.3]{HMTY}).} 
{\it  
These data define a covering group 
$\widetilde{G}$ 
of order {\rm 81}: 
\begin{gather}
\label{2024-01-29-10}
\widetilde{G}=
\big[(Z_{12}\times \Xi_1)\rtimes \Xi_2\big]\rtimes \Xi_3,
\end{gather} 
where  
$Z_{12}:=\langle z_{12}\rangle,\;
\Xi_i:=\langle \xi_i\rangle\;(1\leqslant i\leqslant 3)$. 
With this 
structural datum, $\widetilde{G}$ is identified with 
$G_{81}^{\,10}$ with GAPIdentity\, {\rm [81,10]}. 
 The set of elements of $\widetilde{G}$
is 

$
h=z_{12}^{\;\beta}\,\xi_1^{\,\gamma_1}
\xi_2^{\,\gamma_2}\xi_3^{\,\gamma_3},
\;\,
0{\leqslant} \beta,\gamma_1,\gamma_2,
\gamma_3{\leqslant}2\;({\rm mod}\;3).
$
}

{\bf 2.1. \!Exact identification of $\widetilde{G}$ with $G_{81}^{\,10}$. } 
The above proof of 
identification of $\widetilde{G}$ with $G_{81}^{\,10}$ 
is unsatisfactory for us 
in the sense that we should refer a complete 
list of classification of finite groups with their structural 
data, of order $\leqslant 100$, a result of GAP system using 
 big computers (e.g.\;\cite{WP}). Actually  
 in \cite{Tsur} using 
GAP system, we found a presentation of $G_{81}^{\,10}$, 
different from ours (\ref{2023-11-11-11}) at the point that, 
after translate in our notation,  
the relations\; 
$
\xi_1^{\;3}=1,\;\xi_3^{\;3}=1
$\; 
in (\ref{2023-11-11-11}) is replaced by 
\;$
\xi_1^{\;3}=z_{12},\;\xi_3^{\;3}=z_{12}^{\;\;2}
$.
To fulfill this unsatisfaction, we prove here the following. 

{\bf Theorem 2.2.} {\it Take any pair of integers $a,\,b$ as 
 \,$0\leqslant a,\, b\leqslant 2$, and consider 
\\
\quad{\bf Set of fundamental relations:} 
\begin{gather}
\label{2024-01-23-5}
\left\{
\begin{array}{ll} 
\xi_1^{\;3}=z_{12}^{\:\:a},\;\;\xi_2^{\;3}=1,\;\;
\xi_3^{\;3}=z_{12}^{\;\;b}, 
\\
z_{12}=[\xi_1,\xi_2]\;\;{\rm central},\;\;z_{12}^{\;\;3}=1,
\\
\xi_2=[\xi_1, \xi_3],\quad[\xi_2,\xi_3]=1.
\end{array}
\right.
\end{gather}
Then, with the same {\bf Set of generators} as above, it  
 defines a covering group $G'$ of order {\rm 81}, 
which is naturally isomorphic 
to $\widetilde{G}$.
 }

{\it Proof.}\;  Consider an ideal element 
$\zeta$ which is central, and may  
correspond to $z_{12}^{1/3}$, that is, $\zeta^3=z_{12}$. 
Then, consider an algebraic system $G^\sharp$ 
generated by $\zeta$ 
and $\{\xi_i\;(1\leqslant i\leqslant 3),\, z_{12}\}$ 
under the set of fundamental 
relations ($\ref{2023-11-11-11}'$) below,  
consisting
 (\ref{2023-11-11-11}) and one more line 
\lq\lq\;$\zeta$\;\;central and $\zeta^3=z_{12}$\,.\rq\rq\;\, 
This means 
\\[.5ex]
\quad
{\bf Set of fundamental relations:} 
\\[.2ex]
(\ref{2023-11-11-11}$'$)\hspace*{1.8ex}
$
\mbox{\rm 3 lines in (\ref{2023-11-11-11})\, $\&$ } 
\mbox{\rm \lq\lq\,$\zeta$\;central and $\zeta^3=z_{12}$.\rq\rq}
$
\vskip.5ex

Then, it is proved just like \cite[Theorem 5.3]{HMTY} that 
$G^\sharp$ is a group of order $81\!\times\! 3$, an  
extension of $\widetilde{G}$. 
This bigger group  $G^\sharp$ contains both 
$\widetilde{G}$ and $G'$ 
as its subgroups. 
Taking note on {\bf Set of fundamental relations} (6), 
we consider a map $\Phi$ 
on $G^\sharp$ given by 
$$
\left\{
\begin{array}{ll}
\xi_1\to \xi'_1:=\xi_1\zeta^a,\;\;\xi_3\to 
\xi_3':=\xi_3\,\zeta^b,
\\
\xi_2\to \xi'_2:=\xi_2,\;\; z_{12}\to z_{12},
\;\; \zeta\to \zeta\,.
\end{array}
\right.
$$
Then $\Phi$ is naturally extended 
to an automorphism of $G^\sharp$ and that 
{\bf Set of fundamental relations} (6) for 
$\widetilde{G}$ is mapped 
to $(\ref{2024-01-23-5})$ for $G'$, but with 
{\bf Set of generators} 
$\{\xi'_1,\,\xi'_2,\,\xi'_3,\,z_{12},\,\zeta\}$ instead of. 
\hfill
$\Box$

{\bf Remark 2.3.}\; The above phenomenon shows that, 
under one-step efficient central extensions, 
the resulting extensions are not necessarily  
 unique.

{\bf 2.2. Second step of stairway up to the Sky. } 
2nd step of the present stairway up to the Sky, 
or to 
arrive to
representation group $R(G)$, is to take  
 one-step efficient central extension $H$ 
of $\widetilde{G}$ 
relative to the commuting pair 
$\{\xi_2, \xi_3\}$. 
$H$ is presented as 
\\
\quad
{\bf Set of generators:}\qquad
$\{\eta_1,\eta_2,\eta_3,z_{12},z_{23}\}$\,; 
\\
\quad
{\bf Covering map:} \;$\delta: H \to \widetilde{G}$, 
\begin{gather*}
\label{}
\eta_j\stackrel{\delta}{\rightarrow} \xi_j\;(1\leqslant j\leqslant 3),\;\;
z_{12}\stackrel{\delta}{\rightarrow}z_{12},
\\
[\eta_2,\eta_3]=z_{23}\stackrel{\delta}{\rightarrow}
[\xi_2,\xi_3]=1;
\end{gather*}
\quad
{\bf Set of fundamental relations:}\; 
\begin{gather}
\label{2024-01-29-1}\quad
\left\{
\begin{array}{ll}
\eta_j^{\;3}=1\;(1\leqslant j\leqslant 3),\;&
\varphi(\eta_2)\eta_1=z_{12}^{\;-1}\eta_1,
\\
\varphi(\eta_3)\eta_1=z_{12}\,\eta_1\eta_2^{\;2},\;
&
\varphi(\eta_3)\eta_2=z_{23}^{\;-1}\eta_2;
\end{array}
\right.
\end{gather}
\quad
{\bf Set of elements of $H$\,:}\;(supposed to be $R(G)$)
\begin{gather}
\label{2024-01-29-2}\quad
\left\{
\begin{array}{ll}
\hspace*{3.5ex}
h=z_{12}^{\;\beta_1}z_{23}^{\;\beta_2}\,\eta_1^{\,\gamma_1}
\eta_2^{\,\gamma_2}\eta_3^{\,\gamma_3},
\\
0\leqslant \beta_1,\beta_2,\gamma_1,\gamma_2,\gamma_3
\leqslant 2\;\;({\rm mod}\;3).
\end{array}
\right.
\end{gather}
Actually, under the above product rule, 
the pair of\, {\bf 
Set of generators} and\, 
{\bf Set of fundamental relations}, the 
algebraic system $H$ 
gives a group of order $243$.  
This is a representation group $R(G)$, 
since, with\; 
$
Z_{23}:=\langle z_{23}\rangle,\;  
Y_j:=\langle\eta_j\rangle\;(1\leqslant j\leqslant 3), 
$ 
there holds 
\begin{gather*}
Z(H)=Z_{12}\times Z_{23},\quad
[H,H]
=Z_{12}\times Z_{23}\times Y_2,  
\end{gather*}
which means that the condition (\ref{2024-01-25-11}) is satisfied.

{\bf Theorem 2.4} (cf. \cite[Theorem 6.2]{TYMH}).
 {\it 
Denote by $C_3$ a cyclic group of order  {\rm 3}, 
then 
the group $H=R(G)$ has a structure given as 
\begin{gather}
\label{2024-01-29-4}
R(G)=
\big\{[(Z_{12}\times Y_1)\rtimes Y_2]\times Z_{23}\big\}
\rtimes Y_3\;\,
\\
\nonumber
\hspace*{9.2ex}
\cong\big\{[(C_3\times C_3)\rtimes C_3]\times C_3\big\}\rtimes C_3, 
\end{gather}
and is identified with 
$G_{243}^{\;\;\;3}$ of 
GAPIdentity \,{\rm [243,3]}.
}
\\[-2ex]

{\bf 3. Spin types of projective representations of $G$. }  
A projective representation of $G$ naturally comes from 
 a linear representation $\Pi$ of $R(G)$ (cf. \cite[\S 2]{TYMH}). 
If it is irreducible, 
then we define its spin type as a character $\chi=\chi_{M(G)}$ 
of Schur multiplier 
$M(G)$ in the center of 
$R(G)$, given by $\Pi(z)=\chi(z)I\;(z\in M(G))$. 
 For $G\cong G_{27}^{\;3}$, 
$M(G)=Z_{12}\times Z_{23}$ and its characters are given, 
 with $\omega=\exp(2\pi i/3)$, 
\begin{gather}
\label{2024-02-13-1}
\chi_{\varepsilon, \mu}:=\chi_\varepsilon\cdot \chi_\mu,\;
\chi_\varepsilon(z_{12})=\omega^\varepsilon,\;
\chi_\mu(z_{23})=\omega^\mu,
\end{gather}
for $\varepsilon, \mu=0,1,-1$\;(mod 3). 
To construct explicitly matrix realization of spin IRs and 
to calculate their characters, 
 we discuss 
mainly 
according to their spin types, since there exist 
$9\!=\!3^2$ different spin types in total, 
and simple confusions 
must be avoided. 

The case $\chi_\varepsilon=1,\,\chi_\mu=1$ (or $\varepsilon=\mu=0$) 
is called as 
{\it non-spin}, and the case 
$\chi_\varepsilon\ne 1,\,\chi_\mu\ne 1$ as {\it purely spin}, 
and the remaining cases as {\it partially-spin.}

{\bf 4. Another stairway up to the Sky. }  
We used {\it Stairway up to the Sky} consisting 
of 1st step efficient 
central extension with 
$z_{12}$ and 2nd step with $z_{23}$, to 
 arrive to $R(G)$ in Theorem 2.4 above as
\begin{gather}
\label{2024-02-17-1}
G \;\stackrel{\mbox{\rm\footnotesize with\,$z_{12}$\;}}{\;\longrightarrow\;}
\widetilde{G}\;
\stackrel{\mbox{\rm\footnotesize with\,$z_{23}$\;}}{\;\longrightarrow\;}R(G)\,.
\end{gather}

Studying linear IRs of the covering group 
$\widetilde{G}$, in the 
middle, means 
 that we are dealing with projective 
IRs of spin types coming from 
$\widehat{Z_{12}}$, that is, of spin type 
$\chi_{0,0}$ {\it non-spin}, and of spin type 
$\chi_{\varepsilon,0},\,\varepsilon \not\equiv 0,$ 
{\it partially-spin}. 
To arrive to another type of {\it partially-spin} 
case $\chi_{0,\mu},\,\mu \not\equiv 0,$ 
the shortest 
stairway to take is another choice of order of 
steps of {\it one-step efficient central extensions}, 
that is, first introduce $z_{23}$, and second 
introduce $z_{12}$, or    
\begin{gather}
\label{2024-02-17-2}
G \;\stackrel{\mbox{\rm\footnotesize with\,$z_{23}$\;}}{\;\longrightarrow\;}
\overline{G}\;
\stackrel{\mbox{\rm\footnotesize with\,$z_{12}$\;}}{\;\longrightarrow\;}R(G)\,.
\end{gather}
Spin types of linear IRs of the covering group
$\overline{G}$ comes from characters of 
its central subgroup $Z_{23}$, that is, 
$\chi_{0,0}$ of {\it non-spin} and 
$\chi_{0,\mu},\,\mu\not\equiv 0,$ of 
another {\it partially-spin}. 
Note that we can also 
 arrive to the middle covering group $\overline{G}$ 
if we add to (\ref{2024-02-17-1}) 
one more step after $R(G)$, downwards,  
taking quotient by $z_{12}\to 1$, but 
it's a roundabout stairway.
\!\footnote{
$\overline{G}$ is accidentally 
isomorphic to $\widetilde{G}$ as abstract groups. 
}

{\bf 5. IRs of non-spin type $\chi_{0,0}$. } 
IRs of {\it non-spin} type are nothing but linear 
 IRs of $G$ itself.  
 We apply Mackey's 
induced representation.  
It is dealing with locally compact groups, 
but if applied to finite groups, 
the complexity of set theory disappears. 

{\bf Theorem 5.1} (G.W. Mackey, cf.\;\cite{Mack}).  
{\it Let $G=U\rtimes W$ be 
a finite semidirect product group with $U$ abelian. 
Consider $W$-action\footnote{Actually equivalent 
to consider $G$-action.} 
 on the dual\, $\widehat{U}$ and take a 
complete set of representatives $\{\rho_j;j\in J\}$ of 
$W$-orbits in $\widehat{U}/W.$  Let 
$W(\rho_j)$ be the stationary subgroup of $\rho_j$ 
in $W$.  
Take a complete set of 
 representatives $\{\pi_i\;;\;i\in I_j\}$ 
of IRs of $W(\rho_j)$ and induced representations, 
 for $j\in J,\,i\in I_j$,
\begin{gather*}
\Pi(\rho_j,\pi_i):=
{\rm Ind}^{U\rtimes W}_{U\rtimes W(\rho_j)}
\big(\rho_j\otimes \pi_i\big).
\end{gather*}
 Then the set 
of $\Pi(\rho_j,\pi_i)$'s gives a complete set of 
representatives of IRs of $G=U\rtimes W$.
}

{\bf 5.1. $X_3$\,-\,orbits in $\widehat{U}$ for $U=X_1\times X_2$. } 
Now we return to the case of spin type 
$\chi_{0,0}$ and put 
\begin{gather*}
\label{}
G=U\rtimes W\quad{\rm with}\quad
 U:=X_1\times X_2,\;\;W:=X_3.
\end{gather*}
Then $\widehat{U}=\{\rho_{m,n}\,;\,m,n=0,\,1,-1\}$, where 
\begin{gather}
\label{2024-02-17-21}
\rho_{m,n}(x_1)=\omega^m,\quad 
\rho_{m,n}(x_2)=\omega^n,
\end{gather}
and $w:=x_3$ acts on $U$ as 
\begin{gather*}
\label{2024-02-17-2-1}
\varphi(w)x_1=wx_1w^{-1}=x_2^{\,-1}x_1,\quad
\varphi(w)x_2=x_2.
\end{gather*}

{\bf Lemma 5.2.} {\it {\rm (i)} Action of $w=x_3$ on 
the dual $\widehat{U}$ of\, $U$ is given as 
\;${}^w(\rho_{m,n})=\rho_{m+n,n}$.

{\rm (ii)} $W$-orbits in $\widehat{U}$ are given as 
\begin{gather*}
\widehat{U}/\varphi(W)=
\Big(\bigsqcup_{n=1,-1}\!\!\widehat{U}_{*,n}\Big)
\bigsqcup\Big(\bigsqcup_{m=0,\pm 1}\!\!\big\{\rho_{m,0}\big\}\Big),
\end{gather*}
\\[-2ex]
where\quad 
$\widehat{U}_{*,n}:=\big\{\rho_{m,n}\,
;\,m=0,\,1,-1\big\}$\; 
with representative element $\rho_{0,n}$.

{\rm (iii)} Stationary subgroups in $W=X_3$ 
of each representative 
elements in {\rm (ii)} are given as follows: 
\begin{gather}
\label{2024-02-17-11}
W(\rho_{0,n})=\{1\},\quad W(\rho_{m,0})=W=X_3.
\end{gather}

}

{\bf 5.2. IRs of spin type $\chi_{0,0}$. }  
Applying Mackey's theory, we obtain 
a complete set of IRs of $G=U\rtimes X_3$ as follows:
\begin{gather}
\label{2024-02-20-1}
\Pi_{0,n}:={\rm Ind}^{U\rtimes X_3}_{U\rtimes\{1\}}
\rho_{0,n}\quad(n=1,-1\;{\rm mod}\;3),
\\
\label{2024-02-20-2}
\;\;\,\Pi_{m,0,q}:=\rho_{m,0,q}\qquad
(m,q=0,1,-1 \;{\rm mod}\;3), 
\\
\nonumber
\noindent
{\rm where}
\quad
\hspace*{4ex}
\rho_{m_1,m_2,m_3}(x_i):=\omega^{m_i}\;\;(i=1,2,3).
\hspace{6ex}
\end{gather}

To realize $\Pi_{0,n}$ in matrix form, we take 
$W=X_3$ 
as a section for $U\backslash(U\rtimes W)$. 
Consider the complex function space 
${\mathcal F}(W)$ as $W$-right module, then 
the representation space $V(\Pi_{0,n})$ is 
$V(\rho_{0,n})\otimes{\mathcal F}(W)$ and 
action of $y\in G$ on 
 $f\in V(\Pi_{0,n})$ is given 
by 
\begin{gather}
\label{2024-02-22-11}
\big(\Pi_{0,n}(y)f\big)(x)=f(xy)\quad(x\in G).
\end{gather}
Especially for $\delta$-function $\delta_{x_0}$ 
on a point $x_0$, 
\begin{gather}
\label{2024-02-22-2}
\Pi_{0,n}(y)\delta_{x_0}=\delta_{x_0y^{-1}}.
\end{gather}

Now  $w=x_3$, on $(1, w, w^2)$ is  
$$
1 \stackrel{w}{\longrightarrow} w 
\stackrel{w}{\longrightarrow}w^2 
\stackrel{w}{\longrightarrow} 1
$$
and applying (\ref{2024-02-22-2}) for $y=w$, we have 
the following matrix, as 
the contragredient of the above action, 
\begin{gather}
\label{2024-02-20-3}
\Pi_{0,n}(x_3)= 
{\small 
\begin{bmatrix}
0\;&1\;&0
\\
0\;&0\;&1
\\
1\;&0\;&0
\end{bmatrix}
}
\;(=:J).
\end{gather}

{\bf Theorem 5.3.}\; {\it 
{\rm (i)}\; IRs\, $\Pi_{0,n}$ in 
{\rm (\ref{2024-02-20-1})} 
of\, $G$ are 
3-dimensional, and 
realized, together with $\Pi_{0,n}(x_3)$ 
in {\rm (\ref{2024-02-20-3})}, by  
actions of elements $x_1$ and 
$x_2$ given as 
\begin{gather*}
\label{2024-02-20-4}
\Pi_{0,n}(x_1)={\small 
\begin{bmatrix}
1\! \!&0\! \!&0\\
0\!\!&\omega^{-n} \!\!&0\\
0\!\!&0\!\!&\omega^n
\end{bmatrix} },\;\;
\Pi_{0,n}(x_2)={\small 
\begin{bmatrix}
\omega^n \!\! &0\!\!&0\\
0\!\!&\omega^n\!\!&0\\
0\!\!&0\!\!&\omega^n
\end{bmatrix} }.
\end{gather*}

{\rm (ii)}\; IRs $\Pi_{m,0,q}$ in 
{\rm (\ref{2024-02-20-2})} 
is one-dimensional.
}

{\it Proof.}\; (i)\; 
Action of element $x_1$ on the section $W$ is, 
with $\varphi(w)x_1:=wx_1w^{-1}=x_2^{\,-1}x_1$, 
for $0\leqslant i\leqslant 2$, \;
$w^ix_1=\varphi(w)^ix_1\cdot w^i=(x_2^{\,-i}x_1)\cdot w^i$, and
\begin{gather*}
\label{2024-02-20-5}
(1, w, w^2)\;\stackrel{x_1}{\longrightarrow}
\;(x_1\!\cdot 1,\;x_2^{\,-1}x_1\!\cdot w,\;
x_2x_1\!\cdot w^2).
\end{gather*}
Accordingly set of $\delta$-functions on each of 
$(1,\,w,\,w^2)$ are respectively multiplied by 
scalars given as 
\begin{gather*}
\big(\rho_{0,n}(x_1),\,\rho_{0,n}(x_2^{\,-1}x_1),
\,\rho_{0,n}(x_2x_1)\big).
\end{gather*}
Therefore the matrix $\Pi_{0,n}(x_1)$ is 
diagonal with diagonal elements 
$(1,\,\omega^{-n},\,\omega^n)$. 

Since $x_2$ is a central element of $G$, and so 
$\Pi_{0,n}(x_2)$ is a scalar multiple of $I$. 
\hfill
$\Box$

{\bf 5.3. Characters of spin type $\chi_{0,0}$. }  
The set of elements of $G$ is given as 
\begin{gather*}
h(\beta_1,\beta_2,\beta_3)\!:=x_1^{\,\beta_1}x_2^{\,\beta_2}
x_3^{\,\beta_3},\;
 0\leqslant 
\beta_i\leqslant 2\;{\rm mod}\;3\;(\forall i).
\end{gather*}

{\bf Theorem 5.4.} {\it 
Characters $\chi(h\,|\,\Pi_{0,n})$ 
of IRs $\Pi_{0,n}$ is given as follows:
\begin{gather}
\nonumber
\chi\big(h(\beta_1,\beta_2,\beta_3)\,|\,
\Pi_{0,n}\big)\ne 0\hspace{12ex}
\\
\nonumber
\hspace*{8ex}
\mbox{\rm only if\;\;} \beta_1\equiv 0\;\;
{\rm and}\;\;\beta_3\equiv 0\;\;{\rm mod}\;3,
\\
\nonumber
\label{2024-02-20-11}
\chi\big(h(0,\beta_2,0)\,|\,\Pi_{0,n}\big)=
3\cdot \omega^{\beta_2n}\quad(\beta_2=0,\pm 1).
\end{gather} 
}

{\it Proof.}\; It is easy to see from the matrix form of 
\,$\Pi_{0,n}(x_3)$\, that the
condition\quad  
$\beta_3\equiv 0\;\;{\rm mod}\;\,3$\quad 
is necessary. Then it is sufficient to note that
\begin{gather*} 
\hspace*{9ex}
{\rm tr}\Big(\Pi_{0,n}(x_1^{\,\beta_1})\Big)=0\;\;
{\rm if}\;\;\beta_1=\pm 1.\hspace*{9ex}
\Box
\end{gather*}

{\bf Corollary 5.5.}\; {\it Character 
$\chi(h\,|\,\Pi_{0,n})$ 
of IR $\Pi_{0,n}$ is concentrated on the 
central subgroup $X_2$ of 
$G=(X_1\times X_2)\rtimes X_3$.
 }

{\bf 6. IRs of partially-spin type $\chi_{\varepsilon,0},\,\varepsilon\not\equiv 0$. } 
Recall the structure of $\widetilde{G}$ 
from (\ref{2024-01-29-10}) 
in Theorem 2.1 as
$$
\widetilde{G}=
\big[(Z_{12}\times \Xi_1)\rtimes \Xi_2\big]\rtimes \Xi_3. 
$$
Then linear IRs of $\widetilde{G}$ have their spin types 
from $\chi\in\widehat{Z_{12}}$. If $\chi$ is trivial, it 
 corresponds to {\it non-spin} 
type which is worked out in \S 5.2. So we concentrate here 
on {\it partially-spin} types 
 $\chi_{\varepsilon,0},\,\varepsilon=\pm 1,$  
in the dual of $Z_{12}\times Z_{23}$.\;  
Our study of IRs of this case consists of two-steps 
of semidirect product groups. 

As the first step, we 
take  
\begin{gather}
\label{2024-02-21-1}
U:=U_0\!\rtimes \Xi_2\;\;{\rm with}\;\; U_0\!:=
Z_{12}\!\times \Xi_1\;\;{\rm abelian,}\!\!\!\!
\end{gather}
and as the second step, we take 
\begin{gather}
\label{2024-02-21-2}
\widetilde{G}=U\!\rtimes \Xi_3\;\,{\rm with}\;\,
U=U_0\!\rtimes \Xi_2\;\,\mbox{\rm non-abelian}.\!\!\!\!
\end{gather}

We apply in \S 6.1 Mackey's theory for the 1st 
semidirect product  
(\ref{2024-02-21-1}).  
Then, for the 2nd semidirect 
product (\ref{2024-02-21-2}), we apply in \S 6.2  
the so-called classical method. 
Characters of these IRs of {\it partially-spin} type 
$\chi_{\varepsilon,0},\;\varepsilon \not\equiv 0,$ 
are treated in \S 6.3.

{\bf 6.1. Dual of subgroup $U=U_0\rtimes \Xi_2.$ } 
Action of the group $\Xi_2$ on $U_0$ is given as 
\begin{gather}
\label{2024-01-30-22}
\varphi(\xi_2)z_{12} =z_{12},\quad
\varphi(\xi_2)\xi_1 =z_{12}^{\;-1}\xi_1.
\end{gather}
The dual of 
$Z_{12}$ and that of $\Xi_1$ are given as 
\begin{gather*} 
\label{2024-01-29-12}
\left\{
\begin{array}{rl}
\widehat{Z_{12}}&=\{\chi_\varepsilon\;;\;\varepsilon=0,\pm 
1\},\hspace{4.5ex}\chi_\varepsilon(z_{12}):=\omega^\varepsilon,\;
\\[.4ex]
\widehat{\Xi_1}&=\{\rho_m\;;\;\;m=0,\pm 1\},
\quad\rho_m(\xi_1):=\omega^m,
\end{array}
\right.
\end{gather*}
where $\varepsilon$ and $m$ are counted in mod 3. 
For an IR $\pi$ of $U=U_0\rtimes \Xi_2$, 
we call a character 
$\chi\in \widehat{Z_{12}}$ {\it spin type} of $\pi$ 
if $\pi(z_{12})=\chi(z_{12})I$.

{\bf Lemma 6.1.}\; {\it 
Action of\, $\Xi_2$ on the dual\, 
$\widehat{U_0}=\widehat{Z_{12}}\times \widehat{\Xi_1}$ of 
abelian subgroup $U_0$ is given as follows: put 
$\pi_{\varepsilon,m}:=\chi_\varepsilon\cdot \rho_m$, then  
\begin{gather}
\label{2024-01-29-15}
{}^{\xi_2}\big(\pi_{\varepsilon,m}\big)
=\pi_{\varepsilon,m+\varepsilon}\quad
\mbox{\rm ($m+\varepsilon$, mod\;3)}.
\end{gather}
}

{\it Proof.}\quad ${\displaystyle 
{}^{\xi_2}\big(\pi_{\varepsilon,m}\big)(\xi_1)
=\pi_{\varepsilon,m}\big(\varphi(\xi_2)^{-1}\xi_1\big)
}
\\
\hspace*{13ex}
{\displaystyle 
=\pi_{\varepsilon,m}(z_{12}\xi_1)=\omega^\varepsilon\omega^m
=\omega^{m+\varepsilon}.
}$\hfill
$\Box$

{\bf Lemma 6.2.}\; {\it 
$\Xi_2$-orbits in $\widehat{U_0}$  
are given as  
\begin{gather}
\label{2024-01-30-1}
\nonumber
\widehat{U_0}/\varphi(\Xi_2)= 
\Big(\bigsqcup_{\varepsilon=1,-1}\!\!
\mho_\varepsilon\Big)\bigsqcup\Big(\bigsqcup_{m=0,1,-1}
\!\!
\big\{\pi_{0,m}\big\}\Big),
\end{gather}
\\[-2ex]
where\quad 
$\mho_\varepsilon=
\big\{\pi_{\varepsilon,m}\;;\;m=0,1,-1\big\}$. 
The orbit\, $\mho_\varepsilon$ is represented by 
$\pi_{\varepsilon,0}$, and the stationary subgroup  
 $\Xi_2(\pi_{\varepsilon,0})=\{1\}.$\; 
Other orbits $\{\pi_{0,m}\}$ correspond to 
{\it non-spin case.}
}

Under these Lemmas, we can apply Theorem 5.1 
to obtain explicitly a complete set of representatives 
of IRs of $U=U_0\rtimes \Xi_2$ of spin type $\chi_\varepsilon$.

{\bf Theorem 6.3.} (IRs of $U$ of spin type $\chi_\varepsilon$)
 {\it 
 If an IR $\pi$ of\, $U=U_0\rtimes \Xi_2$ is of spin type  
$\chi_\varepsilon,\,\varepsilon =\pm 1$, then   
$\pi$ is realized by inducing IR  
 $\pi_{\varepsilon,0}$ from $U_0 \times\{1\}$ 
up to $U=U_0\rtimes \Xi_2$ as \;\;
$\begin{displaystyle}  
P_{\varepsilon,0}:=
{\rm Ind}^{U_0\rtimes\Xi_2}_{U_0\times\{1\}}
\pi_{\varepsilon,0}\,,
\end{displaystyle}
$
\begin{gather}
\label{023-08-18-1}
P_{\varepsilon,0}(z_{12})=\omega^{\varepsilon}I,
\quad
\\ 
\label{2024-02-02-11}
P_{\varepsilon,0}(\xi_1){=}{\small 
\begin{bmatrix}
1\;&0\;&0\;
\\
0\;&\omega^{-\varepsilon}&0\;
\\
0\;&0\;&\omega^{\varepsilon}
\end{bmatrix} },
\quad
P_{\varepsilon,0}(\xi_2) = J,
\end{gather}
where $J$ is given in {\rm (\ref{2024-02-20-3})}.
}

{\it Proof.}\;  
Since the stationary subgroup is 
$\Xi_2(\pi_{\varepsilon,0})=\{1\}$, 
we can take $\Xi_2$ for a section of 
  $U_0\backslash U=U_0\backslash(U_0\rtimes\Xi_2)$, as 
$\Xi_2$-right modules, then  
the representation space 
$V(P_{\varepsilon,0})$ 
is isomorphic to the tensor product 
of $V(\pi_{\varepsilon,0})$ 
and the complex 
function space ${\mathcal F}(\Xi_2)$. As a basis  
over 
${\mathbb C}$\, for the latter, 
we take $\delta$-functions on 
each elements of $1,\,\xi_2,\,\xi_2^{\,2}$. 
We can follow similar calculations as those 
in the proof of 
Theorem 5.3 above.
\hfill
$\Box$

Similarly as in Corollary 5.5, we see 

{\bf Lemma 6.4.} {\it Character 
$\chi(h\,|\,P_{\varepsilon,0})$ 
of IR $P_{\varepsilon,0}$ of 
$U=(Z_{12}\times \Xi_1)\rtimes \Xi_2$ 
is concentrated on $Z_{12}$.
 }

{\bf 6.2. IRs of partially-spin type $\chi_{\varepsilon,0}$. }  
For the covering $\widetilde{G}\cong G_{81}^{\;10}$, 
we use its semidirect product structure 
$\widetilde{G}=U\rtimes W$, where  
$W:=\Xi_3=\langle \xi_3\rangle$ acts on 
$U=U_0\rtimes \Xi_2$. 
We 
apply classical method in \cite{THirai}. 
In \cite[\S 3,1]{TYMH}, 
we explained this method of constructing 
a complete set of representatives of IRs of a 
finite semidirect product group $U\rtimes W$ 
rather in detail, tracing from (Step 1) 
up to (Step 4). Please refer loc.\;cit. for  
several notations etc.

{\bf 6.2.1. $\Xi_3$-orbits in the dual of $U=U_0\rtimes \Xi_2.$ }  
\hspace*{3.2ex} {\bf (Step 1).}  Put $w:=\xi_3$, then 
we should check $W$-orbits 
in the dual $\widehat{U}$. Since we are concerned with 
spin type $\chi_{\varepsilon,0}$, it is restricted to 
study only the part of $\widehat{U}$ which have spin 
type $\chi_\varepsilon$ with respect to $Z_{12}$. 
From Theorem 6.3 above, we see that such orbit in 
$\widehat{U}/W$ consists of unique 
element $[\rho],\;\rho:=P_{\varepsilon,0}$. 

{\bf Lemma 6.5.} {\it The stationary subgroup 
$W([\rho]),$  $\rho=P_{\varepsilon,0},$ in $W$ is equal to $W$.
}

{\it Proof.}\; The character $\chi(h\,|\,P_{\varepsilon,0})$ 
is supported by central subgroup $Z_{12}$,  
whence it is invariant under 
$w(u):=\varphi(w)u=wuw^{-1}\;(u\in U)$. \hfill
$\Box$

{\bf (Step 2).}  For $w$ and $w^2 \in W=W([\rho])$, we determine 
 an intertwining operator $J_\rho(w^i)$ given as   
\begin{eqnarray}
\label{2024-02-29-1}
\quad\;\:
\rho\big(w^i(u)\big)J_\rho(w^i)=
J_\rho(w^i)\,\rho(u)\;\,(\forall u\in U), 
\end{eqnarray}
for $0\leqslant i\leqslant 2.$ 
Then each $J_\rho(w^i)$ is determined up to a non-zero 
scalar factor.  
 According to 
fundamental relations (6), we should solve 
the equation (\ref{2024-02-29-1}) using 
(\ref{023-08-18-1})--(\ref{2024-02-02-11}).
Action of $w=\xi_3$ is 
\begin{gather*}
\varphi(\xi_3)\,\xi_1=\xi_2^{\,-1}\xi_1,\;\;
\varphi(\xi_3)^2\xi_1=\xi_2\xi_1,\;\;
\varphi(\xi_3)\,\xi_2=\,\xi_2.
\end{gather*} 


{\bf Lemma 6.6.}\; {\it 
 Put $J$ and $K$ as 
\begin{gather}
\label{2024-02-03-2}
\nonumber
J= {\small 
\begin{bmatrix}
0\;&1\;&0\;
\\
0\;&0\;&1\;
\\
1\;&0\;&0\;
\end{bmatrix}},
\quad
K={\small 
\begin{bmatrix}
0\;&0\;&1\;
\\
1\;&0\;&0\;
\\
0\;&1\;&0\;
\end{bmatrix} }.
\end{gather}
Then, $J^2=K,\,K^2=J,\,J^3=K^3=I$.  
For $\rho=P_{\varepsilon,0}$,  
\,$\rho\big(w^i(\xi_2)\big)=\rho(\xi_2)=J\;(\forall i)$, 
and {\rm (\ref{2024-02-29-1})} 
 for 
$u=\xi_2$ is 
\begin{gather}
\label{2024-03-01-1}
J\cdot J_\rho(w^i)=J_\rho(w^i)\cdot J 
\quad(\forall i).
\end{gather}
Let\, $Y$ be a square matrix 
of order 3, commuting with $J$, 
then $Y=\alpha I+\beta J+\gamma 
K,\,\alpha,\beta,\gamma\in {\mathbb C}$, and 
 ${\rm tr}(Y)=3\alpha,\; 
{\rm det}\,Y=\alpha^3+\beta^3+\gamma^3-3\alpha\beta 
\gamma$. 
 }

\vskip1ex


Put $\rho(\xi_1)
={\rm diag}(1,\omega^{-\varepsilon},\omega^\varepsilon)$, 
then Eq.\;(\ref{2024-02-29-1})  
for $u=\xi_1$ is rewritten as 
\begin{gather*}
\begin{array}{rl}
(*)\qquad\hspace{2ex}
\rho(\xi_1)J_\rho(w)&=J\,J_\rho(w)\rho(\xi_1),\hspace{10ex}
\\
(**)\qquad
\rho(\xi_1)J_\rho(w^2) &=J^2J_\rho(w^2)\rho(\xi_1).
\hspace{10ex}
\end{array}
\end{gather*}

{\bf Lemma 6.7.} {\it Let $J_\rho(w^i)\;(i=0,1,2)$ be a 
projective representation of\, $W$, for which commutativity  
{\rm (\ref{2024-03-01-1})} holds. If $J_\rho(w)$ 
satisfies \mbox{ \rm Eq.\,$(*)$}, then $J_\rho(w)^2$ 
and 
$J_\rho(w)^3$ satisfy respectively the following 
\begin{gather*}
\begin{array}{rl}
(\dag)\qquad\hspace{2ex}
\rho(\xi_1)J_\rho(w)^2&=J^2J_\rho(w)^2\rho(\xi_1),\hspace{10ex}
\\
(\ddag)\qquad\hspace{2ex}
\rho(\xi_1)J_\rho(w)^3 &=J_\rho(w)^3\rho(\xi_1).
\hspace{10ex}
\end{array}
\end{gather*}

}

Thus we see that it is sufficient to solve $(*)$, with 
normalization $J_\rho(w)^3=I$.

{\bf (Step 3).} Here $W([\rho])=W$, and $J_\rho(w)$ 
gives a linear representation of $W$, and the subgroup 
$H:=U\rtimes W([\rho])$ is already equal to 
$\widetilde{G}=U\rtimes W$ itself. So, 
 $\pi^0=\rho\cdot J_\rho$ is an IR of $H=\widetilde{G}$.

{\bf (Step 4).} (cf. \cite[\S 3.1]{TYMH})  
Take an IR $\pi^1$ of $W([\rho])=W\cong C_3$ as $\pi^1=\chi_r,\;
\chi_r(w):=\omega^r,\;
r=0,1,2$, 
and consider it as a representation of $H=\widetilde{G}$ 
through \;$\widetilde{G} \to\, W \cong \widetilde{G}/U.$ 
 Then, inner tensor product 
\begin{gather}
\label{024-03-03-1}
\Pi_{\varepsilon,0;\,r}:=\pi^0\boxdot \pi^1
\end{gather} 
is an IR of $\widetilde{G}=U\rtimes W$ ($\because\; H=\widetilde{G}$).

{\bf Theorem 6.8.} (IRs of spin type $\chi_{\varepsilon,0}$) 
 {\it  
A complete set of representatives of 
IRs of $G\cong G_{27}^{\;\,3}$ of spin type 
$\chi_{\varepsilon,0},\;\varepsilon =1, -1,$ 
is given by 
\begin{gather*}
\label{2024-03-03-11}
\Pi_{\varepsilon,0;\,r}\qquad (r=0,1,2,\;\;{\rm mod\;3}). 
\end{gather*}

}

{\bf 6.3. Characters of spin type $\chi_{\varepsilon,0}$. }  
The set of elements of $\widetilde{G}$ consists of elements 
of the form 
\begin{gather*}
h(\beta,\gamma_1,\gamma_2,\gamma_3)\!:=
z_{12}^{\;\beta}\,\xi_1^{\,\gamma_1}\xi_2^{\,\gamma_2}
\xi_3^{\,\gamma_3},
\end{gather*}
and naturally we need an explicit matrix expression of 
$J_\rho(w),\,w=\xi_3$. We obtain from $(*)$  
 an explicit form of $J_\rho(w)$, and calculate 
spin characters:

{\bf Lemma 6.9.} {\it A solution 
$J_\rho(w)$ 
of equation $(*)$ under 
normalization $J_\rho(w)^3 =I$ is given by 
\begin{gather}
\label{2024-07-18-21}
J_\rho(w)=\alpha
(I{+} \omega^{-\varepsilon} J{+}K)\;\;
\end{gather}
with
$3\alpha^3(1 \!+\! 2\omega^{-\varepsilon})
=1$ or $\alpha \!=\!-i\,\varepsilon/\sqrt{3}$. It is unitary 
and ${\rm det}\big(J_\rho(w)\big)=\omega^\varepsilon.$ 
}

{\bf 7. IRs\;of\;purely-spin\;type $\chi_{\varepsilon,\mu}, 
\varepsilon,\mu\not\equiv\,0$. } 
Representation group $R(G)$ is attained 
from 
 $\widetilde{G}$ through one more step of 
one-step efficient central 
extension $\delta$ as  
in \S 2.2, and 
$$
R(G)=
\big\{\big[(Z_{12}\times Y_1)\rtimes Y_2\big]
\times Z_{23}\big\}\rtimes Y_3. 
$$
Our discussions go almost parallel to those in \S 6, 
 slightly changing appropriately. Put 
\begin{gather}
\label{2024-03-04-12}
\left\{
\begin{array}{ll}
U&\! \! \!:=U_0\!\rtimes Y_2\;\;{\rm with}\;\; U_0\!:=
Z_{12}\!\times Y_1\;{\rm abelian,}\!\!\!\!
\\
\widetilde{U}&\! \! \!:=U\times Z_{23}\;{\rm with}\;\; U\;
{\rm non\!\!-\!\!abelian.}
\end{array}
\right.
\end{gather}

{\bf 7.1. Dual of subgroup $\widetilde{U}$. } 
 Lemmas 6.1 and 6.2 can be translated directly 
by changing $\xi_i\to \eta_i\;(i=1,2)$ and call them 
as {\bf Lemma 7.1} and {\bf Lemma 7.2} respectively. 
Theorem 6.3 is 
transformed for $\widetilde{U}$ as follows. 
Noting that $M(G)=Z_{12}\times Z_{23}$, we call 
a character $\chi_{\varepsilon,\mu}\in \widehat{M(G)}$ 
{\it spin type}. 

{\bf Theorem 7.3.} (IRs of $\widetilde{U}$ 
of spin type $\chi_{\varepsilon,\mu}$)
 {\it \qquad
 If an IR $\pi$ of\, $\widetilde{U}$ is of spin type  
$\chi_{\varepsilon,\mu},\,\varepsilon,\mu \not\equiv 0$, 
then   
$\pi$ is realized as 
\begin{gather}
\label{024-03-04-21}
P_{\varepsilon,\mu}(z_{12})=\omega^\varepsilon I,\;\;
P_{\varepsilon,\mu}(z_{23})=\omega^\mu I,
\quad
\\ 
\label{2024-03-04-22}
P_{\varepsilon,\mu}(\eta_1){=}{\small  
\begin{bmatrix}
1\;&0&0
\\
0\;&\omega^{-\varepsilon}&0
\\
0\;&0&\omega^{\varepsilon}
\end{bmatrix} 
},
\quad
P_{\varepsilon,\mu}(\eta_2)=J. 
\end{gather}
}

{\bf Lemma 7.4.} {\it Character 
$\chi(h\,|\,P_{\varepsilon,\mu})$ 
of IR $P_{\varepsilon,\mu}$ of\, 
$\widetilde{U}$ 
is concentrated on $Z_{12}\times Z_{23}$.
 }

{\bf 7.2. IRs of purely-spin type $\chi_{\varepsilon,\mu}$. }  
For $R(G)\cong G_{243}^{\;\;\,3}$, 
we use its semidirect product structure 
$R(G)=\widetilde{U}\rtimes W$, with
$\widetilde{U}\rtimes W,\;  
W:=Y_3=\langle \eta_3\rangle$.
Here classical method in \cite{THirai} 
is applied. 
In \S 6 before, how to use this method 
 is explained in detail, 
from (Step 1) 
 to (Step 4). 
We can follow it parallelly with 
appropriate changes.

{\bf 7.2.1. $Y_3$-orbits in the dual 
of $\widetilde{U}=U\times Z_{23}$. }   
The difference of situation here from that in \S 6.2.1 
is appearance of central subgroup $Z_{23}$ and 
its character $\chi_\mu(z_{23})=\omega^\mu$. 
Let's add some comments.

In {\bf (Step 1)},  put $w:=\eta_3$. 
$W$-orbits 
in the dual of $\widetilde{U}$ with 
spin type $\chi_{\varepsilon,\mu}$ is 
 $[\rho],\;\rho:=P_{\varepsilon,\mu}$.

\vskip.8ex

{\bf Lemma 7.5.} {\it The stationary subgroup 
$W([\rho]),$  $\rho=P_{\varepsilon,\mu},$ in $W$ 
is equal to $W$.}

 \vskip.8ex

In {\bf (Step 2)}, for $w^i \in W([\rho])=W$,  
 an intertwining operator $J_\rho(w^i)$ is 
determined by Eq.(\ref{2024-02-29-1}). 

In 
{\bf (Step 3)}, $W([\rho])=W$ and $J_\rho(w)$ 
gives a linear representation of $W$, and the subgroup 
$H:=\widetilde{U}\rtimes W([\rho])$ is equal to $R(G)$.
 So, 
 $\pi^0=\rho\cdot J_\rho$ is an IR of $H=R(G)$. 

In 
{\bf (Step 4)}, 
take an IR $\pi^1=\chi_s,\,
\chi_s(w):=\omega^s,\,
s=0,1,2,$ of $W([\rho])=W$. Put  
\begin{gather}
\label{024-03-05-1}
\Pi_{\varepsilon,\mu;\,s}:=\pi^0\boxdot \pi^1.
\end{gather}

{\bf Theorem 7.6.} (IRs of spin type $\chi_{\varepsilon,\mu}$) 
 {\it  
A complete set of representatives of 
IRs of $G\cong G_{27}^{\;\,3}$ of spin type 
$\chi_{\varepsilon,\mu},\;\varepsilon,\mu =1, -1,$ 
is given by 
\begin{gather*}
\label{2024-03-05-11}
\Pi_{\varepsilon,\mu;\,s}\qquad (s=0,1,2,\;\;{\rm mod\;3}). 
\end{gather*}
}

{\bf 8. IRs\;of\;partially-spin\;type
$\chi_{0,\mu},
\mu\not\equiv\,0$. } 
First step 
of {\it Efficient stairway}  
{\it up to the Sky} is one-step 
efficient central extension 
$[\xi_2,\xi_3]=z_{23}$, So, with $\chi_t(w):=\omega^t$, 
similarly as  
in \S 6,  

{\bf Theorem 8.1.} (IRs of spin type $\chi_{0,\mu}$) 
 {\it  
A complete set of representatives of such 
IRs of $G$  
is 
\begin{gather*}
\qquad
\Pi_{0,\mu;\,t}\;\; (t=0,1,2,\;\,{\rm mod\;3})\;\;
\mbox{\rm for}\;\mu\not\equiv 0. 
\end{gather*}
}

\noindent
{\large \bf e-mail addresses:}

\medskip

\hspace{-2mm}Satoe Yamanaka: yamanaka@libe.nara-k.ac.jp

\hspace{-2mm}Tatsuya Tsurii: t3tsuri23@rsch.tuis.ac.jp

\hspace{-2mm}Itsumi Mikami: kojirou@kcn.ne.jp

\hspace{-2mm}Takeshi Hirai: hirai.takeshi.24e@st.kyoto-u.ac.jp

\end{multicols*}

\end{document}